\documentclass[a4paper,10pt]{article}
\usepackage{amsmath,amsthm,amssymb,amscd,epsfig,esint, mathrsfs,graphics,color}
 \usepackage{wrapfig}
\usepackage{verbatim}
\usepackage[english]{babel}
cm \oddsidemargin=.2true cm \topmargin=-.5truecm

\newcommand{\R}{\numberset{R}}

\newcommand{\e}{\varepsilon}
\newcommand{\om}{\Omega}
\newcommand{\blr}{{B_{\lambda r}}}
\newcommand{\bt}{{B_t}}

\theoremstyle{plain}
\newtheorem{thm}{Theorem}[section]

\newtheorem{proposition}[thm]{Proposition}
\newtheorem{lemma}[thm]{Lemma}

\theoremstyle{definition}
\newtheorem{definition}[thm]{Definition}

\def\XXint#1#2#3{{\setbox0=\hbox{$#1{#2#3}{\int}$}
        \vcenter{\hbox{$#2#3$}}\kern-.5\wd0}}

\def\dx{{\mathrm d}x}
\def\d\theta{{\mathrm d}\theta}

\def\eps{\varepsilon}

\def\R{\mathbb{R}}
\def\e{\varepsilon}

\numberwithin{equation}{section} \makeatletter

\renewcommand{\p@enumi}{\thesection.}

\title{\sc{Regularity results for solutions to a class of obstacle problems}
\footnotetext{\hspace{-0.35cm} 2010 \emph{Mathematics Subject
Classification}. 35J87, 49J40.
\endgraf
{\it Key words and phrases}.  }}
\author{ Andrea Gentile, Raffaella Giova, Andrea Torricelli}

\title{\textbf{Regularity results for bounded solutions to obstacle problems with non-standard growth conditions }}
\begin{document}
\maketitle

%
%

\begin{abstract}
In this paper we consider a
class of obstacle problems of the type


\begin{equation*}
\min \left\{\int_{\Omega}f(x, Dv)\, \dx\,:\, v\in
\mathcal{K}_\psi(\Omega)\right\}
\end{equation*}
where $\psi$ is the obstacle,
$\mathcal{K}_\psi(\Omega)=\{v\in u_0+W^{1, p}_{0}(\Omega, \R): v\ge\psi \text{ a.e. in }\Omega\}$, with $u_0 \in W^{1,p}(\Omega)$ a fixed boundary datum, the class of the admissible functions
and the integrand $f(x, Dv)$
satisfies non standard $(p,q)$-growth conditions. \\
We prove higher differentiability results for bounded solutions of
the obstacle problem under dimension-free conditions on the gap between the growth and the ellipticity exponents.
Moreover, also the Sobolev assumption on the partial map
$x\mapsto A(x, \xi)$ is independent of the dimension $n$ and this, in some cases, allows us to manage coefficients in a Sobolev class below the
critical one $W^{1,n}$.

\end{abstract}


\bigskip

\bigskip

\section{Introduction}

We prove higher differentiability results for solutions
to variational obstacle problems of the form

\begin{equation}\label{functionalobstacle}
\min \left\{\int_{\Omega}f(x, Dv)\, \dx\,:\, v\in
\mathcal{K}_\psi(\Omega)\right\},
\end{equation}
where $\Omega$ is a bounded open set of $\R^n$, $n>2$,  $\psi: \Omega
\mapsto [-\infty, +\infty)$ belonging to the Sobolev class $ W^{1,
p}_{\mathrm{loc}}(\Omega)$ is the \emph{obstacle} and

$$\mathcal{K}_\psi(\Omega)=\{v\in u_0+W^{1, p}_{0}(\Omega, \R): v\ge\psi \text{ a.e. in }\Omega\}$$
is the class of the admissible functions, with $u_0 \in W^{1,p}(\Omega)$ a fixed boundary datum.

\bigskip

We shall consider integrands $f$ such that $\xi \mapsto f(x, \xi)$ is $\mathcal{C}^2$ and  there
exists $\tilde{f}:\Omega\times [0, \infty) \to [0, \infty)$ such that

$$f(x,\xi)=\tilde{f}(x, |\xi|).$$

\noindent Moreover, we assume that there exist positive constants $\tilde{\nu}, \tilde{L}$,
exponents $p,q$ with $2\le p< q < p+1< +\infty$ and a parameter $0\le \mu\le 1$  such that the following assumptions are satisfied


$$\langle D_{\xi\xi}f(x,\xi)\lambda,\lambda\rangle\geqslant\tilde{\nu}(\mu^2+|\xi|^2)^{\frac{p-2}{2}}|\lambda|^{2} \eqno{\rm (F1)}$$

 $$| D_{\xi\xi}f(x,\xi)|\leqslant \tilde{L}(\mu^2+|\xi|^2)^{\frac{q-2}{2}}\eqno{\rm (F2)} $$
 for almost every $x\in \Omega$ and every $\xi, \lambda \in \R^n$.\\

\noindent Note that, as proved in \cite{CMMP}, the assumptions (F1) and (F2) and the dependence on the modulus imply that there exists a positive constant $\tilde{\ell}$ such that

 $$\frac{1}{\tilde{\ell}} (|\xi|^2-\mu^2)^{\frac{p}{2}}\le f(x,\xi)\le \tilde{\ell}(\mu^2+|\xi|^2)^{\frac{q}{2}} \eqno{\rm (F3)}$$
 for almost every $x\in \Omega$ and every $\xi \in \R^n$, i.e. the functional $f$ has non-standard growth conditions of $(p,q)$-type as defined and introduced by Marcellini in \cite{Marcellini1, Marcellini2, Marcellini3}
and then widely investigated (see for example \cite{CKP, ELM1,  ELM2}).
\bigskip
\noindent Concerning the dependence on the $x-$variable, we assume that there exists  a non negative function $k(x)\in L^{\frac{p+2}{p-q+1}}$ such that
 $$|D_{x\xi}f(x,\xi)|\leqslant k(x)(\mu^2+|\xi|^2)^{\frac{q-1}{2}}\eqno{\rm (F4)}$$
for almost every $x\in \Omega$ and every $\xi \in \R^n$.

\bigskip

Let us observe that, in case of standard growth conditions, $u\in W^{1,p}_{\mathrm{loc}}(\Omega)$ is a
solution to the obstacle problem \eqref{functionalobstacle} in
$\mathcal{K}_\psi(\Omega)$ if and only if $u \in
\mathcal{K}_\psi(\Omega)$ and $u$ is a solution to the variational
inequality

\begin{equation}\label{variationalinequality}
\int_{\Omega}\left<A(x, Du(x)), D(\varphi(x)-u(x))\right>\dx\ge0\qquad\forall
\varphi\in W^{1,\infty}_{\mathrm{loc}}(\Omega)\, \mathrm{and} \, \varphi \ge \psi,
\end{equation}
where the operator $A(x, \xi): \Omega\times\R^n\to\R^n$ is defined
as follows

\begin{equation*}
A(x, \xi)=D_{\xi}f(x, \xi).
\end{equation*}

\noindent It is clear that, in case of standard growth, a density argument shows the validity of \eqref{variationalinequality} for every $ \varphi \in \mathcal{K}_\psi(\Omega)$.
Here, dealing with non standard growth, it is worth observing that \eqref{variationalinequality} holds true also for solutions to \eqref{functionalobstacle}.
More precisely, due to our assumptions $q-p<1$ on the gap between the ellipticity exponent $p$ and the growth exponent $q$, the
validity of \eqref{variationalinequality} can be easily checked as done at the beginning of the proof of the Theorem \ref{thm1} below.\\
We want to stress that this is not obvious in case of non standard growth conditions: already for unconstrained problems, the relation between
minima and extremals, i.e. solutions of the corresponding Euler Lagrange system, is an issue that requires a
careful investigation (see for example \cite{CKP2, CKP3} and for constrained problems see the very recent paper \cite{EP3}).

\bigskip
From assumptions
$\rm (F1)-\rm (F4)$, we deduce the existence of positive constants $\nu, L, \ell$ such that the following $p$-ellipticity
and $q$-growth conditions are satisfied by the map $A$:

$$\left<A(x, \xi)-A(x, \eta), \xi-\eta\right>\ge\nu|\xi-\eta|^2\left(\mu^2+|\xi|^2+|\eta|^2\right)^\frac{p-2}{2} \eqno{\rm (A1)}$$

$$\left|A(x, \xi)-A(x, \eta)\right|\le
L|\xi-\eta|\left(\mu^2+|\xi|^2+|\eta|^2\right)^\frac{q-2}{2}
\eqno{\rm (A2)}$$

$$
\left|A(x, \xi)\right|\le
\ell\left(\mu^2+|\xi|^2\right)^\frac{q-1}{2}, \eqno{\rm (A3)}$$

$$|D_{x}A(x,\xi)|\leqslant k(x)(\mu^2+|\xi|^2)^{\frac{q-1}{2}}\eqno{\rm (\tilde{A}4)}$$

 for almost every $x\in\Omega$ and for every $\xi, \eta\in\R^n$.\\

Thanks to a characterization of the Sobolev spaces due to Hajlasz \cite{Hajlasz},
we deduce from $\rm (\tilde{A}4)$ that there exists a non-negative function $\kappa\in
L^{\frac{p+2}{p-q+1}}_{\mathrm{loc}}(\Omega)$ such that

$$
\left|A(x, \xi)-A(y, \xi)\right|\le
\left(\kappa(x)+\kappa(y)\right)|x-y|\left(\mu^2+|\xi|^2\right)^\frac{q-1}{2}
\eqno{\rm (A4)}$$

\noindent for almost every $x, y\in\Omega$ and for all $\xi\in \R^n.$
As far as we know, regularity results concerning local minimizers
of integral functionals of the Calculus of Variations under an assumption on the dependence on the $x$-variable of this type, have been
obtained, for the first time, in \cite{KristensenMingione, KuusiMingione}.

\bigskip
The study of the regularity properties of solutions to obstacle problems has been object of intense interest in the last years and it has been usually observed that the regularity of the obstacle influences the regularity of the solutions to the problem: for linear problems the solutions are as regular as the obstacle;
this is no longer the case in the nonlinear setting for general
integrands without any specific structure. Hence along the years, there has been an intense research activity in which
extra regularity has been imposed on the obstacle to balance the
nonlinearity. \\
Here, as we already said, we are interested in higher differentiability results since in case of non standard growth,  many questions
are still open.
In \cite{ CGG, CDF, EP, Gavioli1,  Gavioli2, Gentile3, GI, MaZ, ZZ} the authors analyzed how
an extra differentiability of integer or fractional order of the
gradient of the obstacle provides an extra differentiability  to the gradient of the solutions, also in case of standard growth. However, since no extra differentiability properties for the
solutions can be expected even if the obstacle $\psi$ is smooth,
unless some assumption is given on the $x$-dependence of the
operator $A$, the higher differentiability results for the solutions of systems or for the minimizers of functionals in the case of unconstrained problems (see \cite{BDW,  CUMR, CGGP, Gentile1, Gentile2,  Giova1, Giova2, GP, APdN1, APdN2})
have been useful and source of inspiration also for the constrained case.

It is well known that, for unconstrained problems with $(p,q)$-growth, the boundedness of the minimizers can play a crucial role in order to get regularity for the gradient, under weaker assumptions
on the gap between $p$ and $q$ and on the data of the problem (see \cite{CKP}). Here, we will prove that the same phenomenon happens for the bounded solutions to obstacle problems with $(p,q)$-growth.

\bigskip
More precisely, we prove the following

\begin{thm}\label{thm1}
 Let $u\in \mathcal{K}_\psi(\Omega)$ be a solution to the obstacle problem \eqref{functionalobstacle} and
    let $A(x, \xi)$ satisfy the assumptions (A1)--(A4) with $2\le p< q<\min \{p+1, p^*=\frac{np}{n-p} \}$.
    Then, if $\psi \in L^{\infty}_\mathrm{loc}(\Omega)$ the
    following implication holds

$$D\psi\in W^{1, \frac{p+2}{p+2-q}}_{\mathrm{loc}}(\Omega) \,\, \implies   \,\, \left(\mu^2+\left|Du\right|^2\right)^\frac{p-2}{4}Du\in W^{1,2}_\mathrm{loc}(\Omega),$$
with the following estimate

\begin{align}\label{estimate}
\int_{B_{\frac{R}{4}}}&\left|DV_p(Du(x))\right|^2 \dx
\le \frac{c(\Arrowvert \psi\Arrowvert_{L^{\infty}}^2+\Arrowvert
    u\Arrowvert_{L^{p^*}\left(B_{R}\right)}^2)}{R^\frac{p+2}{2}}\notag\\
& \qquad \qquad \cdot
\int_{B_{R}}\left[1+\left|D^2\psi(x)\right|^\frac{p+2}{p+2-q}+
\left|D\psi(x)\right|^\frac{p+2}{p+2-q}+\kappa^{\frac{p+2}{p-q+1}}+\left|Du(x)\right|^p
\right] \dx.
\end{align}

\end{thm}

We first observe that the assumption of boundedness of the
obstacle $\psi$ is needed to get the boundedness of the solutions (see Theorem \ref{boundedness}). Therefore, if we want to remove the hypothesis $\psi \in L^{\infty}$, it is sufficient to deal with a priori bounded
minimizers. In this case, we can remove also the hypothesis $q< p^*$.\\
Let us compare, now, our result with the previous ones. All previous higher regularity results for solutions to obstacle problem in case of non-standard growth have been obtained under
a Sobolev assumption $W^{1,r}(\Omega)$ with $r\ge n$  on the dependence on $x$ of the operator $A$. Dealing with bounded solutions, we are able to prove our result assuming that the partial map $x \mapsto A(x; \xi)$ belongs to a Sobolev class that is not related to the dimension $n$ but to the ellipticity and the growth exponents $p$ and $q$ of the
functional and this assumption in case $\frac{p+2}{p-q+1}<n$ (i.e. $p<n-2$ and $q< \frac{n-1}{n}p+\frac{n-2}{n}$) improves the higher differentiability result obtained in \cite{Gavioli2}. Moreover, our result is obtained under a weaker assumption also on the gradient of the obstacle, indeed previous result assumed $\psi \in W^{1, 2q-p}$ (see \cite{Gavioli2}) while our hypothesis is $\psi \in W^{1, \frac{p+2}{p+2-q}}$, and under our assumption on the gap, i.e. $q-p<1$, it results $W^{1, 2q-p}\hookrightarrow W^{1, \frac{p+2}{p+2-q}}$.\\
Note that for $p=q$ we recover exactly our previous result (\cite{CGG}) concerning the obstacle problem with standard growth. \\
On the other hand, our result extends to the solutions of constrained problems the higher differentiability result obtained in \cite{CGGP} for the solutions to unconstrained problems in case of the integrand $f$
is uniformly convex only at infinity.

\bigskip
 \noindent In order to prove Theorem
\ref{thm1}, we first verify the validity of the variational inequality also in the case of non standard growth and then we  combine an a priori estimate for
the second derivatives of the local solutions, obtained using the
difference quotient method, with a suitable approximation argument.
The local boundedness of the obstacle, and then of the solutions, allows us to use two interpolation inequalities
that give the higher local integrability $L^{\frac{2(p+2)}{p+2-q}}$ for the gradient
of the obstacle and the higher local integrability  $L^{p+2}$ of the gradient
of the solutions. Such higher integrability is the key tool in order
to weaken the assumption on $\kappa$ that is the function that control the dependence on $x$-variable of the operator $A$ .\\
We conclude observing that, if the minimizer $u$ is assumed a priori in a Lebesgue space $L^r$ with $r> \frac{np}{n-p-2}$ instead  of assuming $u \in L^{\infty}$
 the interpolation
inequality of Lemma \ref{lemma5GP} still gives a higher integrability result for $Du$, i.e. $Du \in L^{\frac{r}{r+2}(p+2)}$. Such higher integrability allows us to obtain the same higher
differentiability result of Theorem \ref{thm1} assuming $\kappa \in L^{\frac{r}{(r-p)}\frac{(p+2)}{p-q+1}} $. We’d like to point out that for $p<n-2$ and $q< \frac 1 n (n-\frac{r}{r-p})p+ \frac 1 n (n-2\frac{r}{r-p})$ we get
$\frac{r}{(r-p)}\frac{(p+2)}{p-q+1} < n$ that means that we obtain the regularity result again under a Sobolev assumption on the dependence on the $x$-variable  below the
critical one $W^{1,n}$.

\bigskip

\section{Notations and preliminary results}

In this paper we shall denote by $C$ or $c$ a general constant that may vary on different occasions, even within the same
line of estimates. Relevant dependencies on parameters and special
constants will be suitably emphasized using parentheses or
subscripts. With the symbol $B(x,r)=B_r(x)=\{y\in
\R^n:\,\, |y-x|<r\}$ we will denote the ball centered at $x$ of
radius $r$. We
shall omit the dependence on the center  when no confusion arises.

\bigskip

Here we recall some results that will be useful in the following.\\
The main tools in the proof of Theorem \ref{thm1} are the following Gagliardo-Nirenberg-type inequalities that we
state as  lemmas. 
The proof of inequality \eqref{2.1GP} can be found in \cite[Appendix A]{CKP} while inequality \eqref{2.2GP} is
a particular case ($p(x) \equiv p$,  for all $x$) of \cite[Lemma 3.5]{GiannettiPassa}. For the proof of \eqref{brezis} see for example \cite{Nirenberg}.

\begin{lemma}\label{lemma5GP}
    For any $\phi\in C_0^1(\Omega)$ with $\phi\ge0$, and any $C^2$ map $v:\Omega\to\R^N$, we have

    \begin{align}\label{2.1GP}
    \int_\Omega&\phi^{\frac{m}{m+1}(p+2)}(x)|Dv(x)|^{\frac{m}{m+1}(p+2)}\dx\notag\\
    \le&(p+2)^2\left(\int_\Omega\phi^{\frac{m}{m+1}(p+2)}(x)|v(x)|^{2m}\dx\right)^\frac{1}{m+1}\cdot\left[\left(\int_\Omega\phi^{\frac{m}{m+1}(p+2)}(x)\left|D\phi(x)\right|^2\left|Dv(x)\right|^p \dx\right)^\frac{m}{m+1}\right.\notag\\
    &\left.+n\left(\int_\Omega\phi^{\frac{m}{m+1}(p+2)}(x)\left|Dv(x)\right|^{p-2}\left|D^2v(x)\right|^2\dx\right)^\frac{m}{m+1}\right],
    \end{align}

    for any $p\in(1, \infty)$ and $m>1$. Moreover, for any $\mu\in[0,1]$

    \begin{align}\label{2.2GP}
        \int_{\Omega}&\phi^2(x)\left(\mu^2+\left|Dv(x)\right|^2\right)^\frac{p}{2}\left|Dv(x)\right|^2 \dx\notag\\
        \le&c\Arrowvert v\Arrowvert_{L^\infty\left(\mathrm{supp}(\phi)\right)}^2\int_\Omega\phi^2(x)\left(\mu^2+\left|Dv(x)\right|^2\right)^\frac{p-2}{2}\left|D^2v(x)\right|^2 \dx\notag\\
        &+c\Arrowvert v\Arrowvert_{L^\infty\left(\mathrm{supp}(\phi)\right)}^2\int_\Omega\left(\phi^2(x)+\left|D\phi(x)\right|^2\right)\left(\mu^2+\left|Dv(x)\right|^2\right)^\frac{p}{2}\dx,
    \end{align}

    for a constant $c=c(p).$
    \end{lemma}

   \begin{lemma}\label{lemmaBrezis}
    Let $u \in L^p(\Omega) \cap W^{2,r}(\Omega)$ with $1 \le p\le \infty$ and $1 \le r \le \infty$. Then $ u \in W^{1,q}(\Omega)$ where $q$ is such that
    $\frac{1}{q}= \frac{1}{2} \left( \frac{1}{p} + \frac{1}{r}\right) $ and

    \begin{equation}\label{brezis}
      \|  Du \|_{L^q} \le C  \|  u \|_{W^{2,r}}^{\frac 1 2}  \|  u \|_{L^p}^{\frac 1 2 }
    \end{equation}

    \end{lemma}

\bigskip
The following is an higher differentiability result to the solutions to \eqref{functionalobstacle} when the energy density function $f$ satisfies standard growth conditions. The proof can be found in \cite{CGG}.
\begin{thm}\label{thmCGG}
    Let $A(x, \xi)$ satisfy the conditions (A1)--(A4) with $p=q\ge2$ and
    let $u\in \mathcal{K}_\psi(\Omega)$ be a solution to the obstacle problem \eqref{variationalinequality}.
    Then, if $\psi \in L^{\infty}_\mathrm{loc}(\Omega)$ the
    following implication

$$D\psi\in W^{1, \frac{p+2}{2}}_{\mathrm{loc}}(\Omega) \,\, \implies   \,\, \left(\mu^2+\left|Du\right|^2\right)^\frac{p-2}{4}Du\in W^{1,2}_\mathrm{loc}(\Omega),$$
 holds true.
 \end{thm}
\noindent Next result has been proved in \cite[ Theorem 1.1]{CEP}

\begin{thm}\label{boundedness}
    Let u in $K_{\psi}(\Omega)$ be a solution of
    \eqref{functionalobstacle} under the assumptions  \rm(A1) and \rm (A2) with $2\le p \le q$ such that
    $$
\begin{array}{llll}
       p \le q < p^\ast=\frac{np}{n-p} & \rm if & p<n
         \\
         \\
        p \le q < \infty  & \rm if & p \ge n \;\;
         \end{array}
$$
If the obstacle $\psi \in L^{\infty}_{\mathrm{loc}}(\Omega)$, then $u
    \in L^{\infty}_{\mathrm{loc}}(\Omega)$ and the following estimate
    \begin{equation}\label{est boundedness}
    \sup_{B_{R/2}} |u| \le \left[ \sup_{B_R} |\psi|+ \left(
    \int_{B_R}|u(x)|^{p^*} \dx\right)\right]^{\gamma}
    \end{equation}
    holds for every ball $B_R \Subset \Omega$, for $\gamma(n,p,q)>0$ and
    $c=c(\ell, \nu, p, q, n)$.
\end{thm}

\bigskip

\noindent We will use the auxiliary function $V_p:\R^n\to\R^n$,
defined as

\begin{equation}\label{Vp}
V_p(\xi):=\left(\mu^2+|\xi|^2\right)^\frac{p-2}{4}\xi,
\end{equation}

\noindent for which the following estimates hold (see \cite{ Giusti} ).

\begin{lemma}\label{lemma6GP}
Let $1<p<\infty$. There is a constant $c=c(n, p)>0$ such that

\begin{equation}\label{lemma6GPestimate1}
    c^{-1}\left(\mu^2+|\xi|^2+|\eta|^2\right)^\frac{p-2}{2}\le\frac{\left|V
_p(\xi)-V_p(\eta)\right|^2}{|\xi-\eta|^2}\le c\left(\mu^2+|\xi|^2+|\eta|^2\right)^\frac{p-2}{2},
    \end{equation}

\noindent for any $\xi, \eta\in\R^n$ and $\xi \ne \eta$.
Moreover, for a $C^2$ function $g$, there is a constant $C(p)$ such that

\begin{equation}\label{lemma6GPestimate2}
C^{-1}\left|D^2g\right|^2\left(\mu^2+\left|Dg\right|^2\right)^\frac{p-2}{2}\le\left|D\left(V_p(Dg)\right)\right|^2\le
C\left|D^2g\right|^2\left(\mu^2+\left|Dg\right|^2\right)^\frac{p-2}{2}.
\end{equation}
\end{lemma}

\bigskip

Now we state a well-known iteration lemma (the proof can be found for example in
\cite[Lemma 6.1]{Giusti}).

\begin{lemma}[Iteration Lemma]\label{iteration}
    Let $h: [\rho, R]\to \R$ be a nonnegative bounded function, $0<\theta<1$, $A, B\ge0$ and $\gamma>0$. Assume that

    $$
    h(r)\le\theta h(s)+\frac{A}{(s-r)^\gamma}+B
    $$

    for all $\rho\le r<s\le R_0<R.$ Then

    $$
    h(\rho)\le \frac{cA}{(R_0-\rho)^\gamma}+cB,
    $$

    where $c=c(\theta, \gamma)>0$.
\end{lemma}

\subsection{Difference quotient}
\medskip
\noindent \noindent In order to get the regularity of the solutions
of the problem \eqref{functionalobstacle}, we shall use the
difference quotient method. We recall here the definition
and basic results.
\begin{definition}
Given $h\in\mathbb{R}^n$, for every function
$F:\mathbb{R}^{n}\to\mathbb{R}$ the finite difference operator is
defined by
$$
\tau_{h}F(x)=F(x+h)-F(x).
$$
\end{definition}

\par

We recall some properties of the   finite difference operator that
will be needed in the sequel. We start with the description of some
elementary properties that can be found, for example, in
\cite{Giusti}.

\bigskip

\begin{proposition}\label{findiffpr}

Let $F$ and $G$ be two functions such that $F, G\in
W^{1,p}(\Omega)$, with $p\geq 1$, and let us consider the set
$$
\Omega_{|h|}:=\left\{x\in \Omega : dist(x,
\partial\Omega)>|h|\right\}.
$$
Then
\begin{itemize}
\item[$(d1)$] $\tau_{h}F\in W^{1,p}(\Omega_{|h|})$ and
$$
D_{i} (\tau_{h}F)=\tau_{h}(D_{i}F).
$$
\item[$(d2)$] If at least one of the functions $F$ or $G$ has support contained
in $\Omega_{|h|}$ then
$$
\int_{\Omega} F(x)\, \tau_{h} G(x)\, \dx =\int_{\Omega} G(x)\, \tau_{-h}F(x)\,
\dx.
$$
\item[$(d3)$] We have
$$
\tau_{h}(F G)(x)=F(x+h )\tau_{h}G(x)+G(x)\tau_{h}F(x).
$$
\end{itemize}
\end{proposition}

\noindent The next result about finite difference operator is a kind
of integral version of Lagrange Theorem.
\begin{lemma}\label{le1} If $0<\rho<R$, $|h|<\frac{R-\rho}{2}$, $1 < p <+\infty$,
 and $F, DF\in L^{p}(B_{R})$ then
$$
\int_{B_{\rho}} |\tau_{h} F(x)|^{p}\ dx\leq c(n,p)|h|^{p}
\int_{B_{R}} |D F(x)|^{p}\ \dx .
$$
Moreover
$$
\int_{B_{\rho}} |F(x+h )|^{p}\ dx\leq  \int_{B_{R}} |F(x)|^{p}\ \dx .
$$
\end{lemma}

\noindent We conclude this section recalling this result that is
proved in \cite{Giusti}.

\begin{lemma}\label{Giusti8.2}
    Let $F:\R^n\to\R^N$, $F\in L^p(B_R)$ with $1<p<+\infty$. Suppose that there exist $\rho\in(0, R)$ and $M>0$ such that

    $$
    \sum_{s=1}^{n}\int_{B_\rho}|\tau_{s, h}F(x)|^p\dx\le M^p|h|^p
    $$

    for every $h<\frac{R-\rho}{2}$. Then $F\in W^{1,p}(B_R, \R^N)$. Moreover

    $$
    \Arrowvert DF \Arrowvert_{L^p(B_\rho)}\le M.
    $$
\end{lemma}

\bigskip

\subsection{Approximation Lemma}
\medskip
We report a Lemma which will be the main tool in the second part of the proof of our main result. For the proof of this Lemma
we refer to \cite{CGM}.

\begin{lemma}\label{approximation}
Let $f:\Omega\times \R^n \to [0, \infty)$ be a Carathéodory function such that $\xi \mapsto f(x, \xi)$ is $\mathcal{C}^2$  and
there exists $\tilde{f}:\Omega\times [0, \infty) \to [0, \infty)$ such that $f(x,\xi)=\tilde{f}(x, |\xi|)$. Moreover, let us assume that $f$ satisfies assumptions (F1)–(F4).
Then there exists a sequence $(f_\e)_{\e}$ of Carathéodory functions $f_\e:\Omega\times \R^n \to [0, \infty)$, monotonically convergent to $f$, such that

\begin{itemize}
  \item [\rm{(i)}] for a.e. $x \in \Omega$, for every $\xi\in \R^n$ and for every $\e_1 < \e_2$, we have
\[
f_{\e_2}(x, \xi)\le f_{\e_1}(x, \xi)\le f(x, \xi)
\]

  \item [\rm{(ii)}] there exists $\bar{\nu}>0$ depending only on $p$ and $\tilde{\nu}$ such that

 $$ \langle D_{\xi\xi}f_\e(x,\xi)\lambda,\lambda\rangle\geqslant\bar{\nu}(\mu^2+|\xi|^2)^{\frac{p-2}{2}}|\lambda|^{2} $$

  for a.e. $x \in \Omega$, for every $\xi\in \R^n$,

  \item [\rm{(iii)}] there exist $K_0, K_1$ independent of $\e$ and $\bar{K}_1$ depending on $\e$ such that

   $$K_0 (|\xi|^p-\mu^2)\le f_\e(x,\xi)\le K_1(\mu^2+|\xi|^q),$$
   $$f_\e(x,\xi)\le \bar{K}_1(\e)(\mu^2+|\xi|^p),$$
for a.e. $x \in \Omega$, for every $\xi\in \R^n$,

 \item [\rm{(iv)}]  there exists a constant $C(\e)>0$ such that
  $$|D_{x\xi}f_\e(x,\xi)|\leqslant k(x)(\mu^2+|\xi|^2)^{\frac{q-1}{2}}$$
   $$|D_{x\xi}f(x,\xi)|\leqslant C(\e)k(x)(\mu^2+|\xi|^2)^{\frac{p-1}{2}}$$

   for a.e. $x \in \Omega$, for every $\xi\in \R^n$.
\end{itemize}

\end{lemma}

\section{Proof of the Theorem \ref{thm1}}

The proof of the theorem is obtained in two steps: first we establish the a priori estimate and then
we conclude through an approximation argument.
\begin{proof}

{\bf Step 1: The a priori estimate.}

In order to get the a priori estimate we first need to prove  the validity of the variational inequality \eqref{variationalinequality} also in the case of non-standard growth conditions.

\noindent Suppose that $u$ is a local
solution to the obstacle problem in $\mathcal{K}_\psi(\Omega)$ such
that
\begin{equation}\label{apriori_estimate_second_derivatives}Du \in W^{1,2}_\mathrm{loc}(\Omega)\text{\qquad and\qquad}\left
(\mu^2+\left|Du\right|^2\right)^\frac{p-2}{4}Du\in
W^{1,2}_\mathrm{loc}(\Omega).
\end{equation}

\noindent Thanks to our assumptions on the exponents $p$ and $q$ we can deduce from Theorem \ref{boundedness} that the solution $u$ to \eqref{functionalobstacle} is bounded. Such boundedness, with the a priori assumption \eqref{apriori_estimate_second_derivatives}  on the second derivatives of $u$, allows us to apply  Lemma \ref{lemma5GP} to get the higher integrability   $Du\in
L^{p+2}_\mathrm{loc}(\Omega)$.\\
Concerning the obstacle $\psi$, by the assumptions $\psi \in L^{\infty}(\Omega)$ and $D^2\psi \in L^{\frac{p+2}{p+2-q}}(\Omega)$, applying Lemma \ref{lemmaBrezis}, we have $D\psi \in L^{\frac{2(p+2)}{p+2-q}}(\Omega)\hookrightarrow L^{p+2}(\Omega)$.

\bigskip

\noindent Note that
$Du\in L^{p+2}_\mathrm{loc}(\Omega)$ (and then, obviously, $u \in W^{1, q}_\mathrm{loc}(\Omega)$)  implies that the variational
inequality \eqref{variationalinequality}, by a simple density
argument, holds
true for every $\varphi \in W^{1, q}_\mathrm{loc}(\Omega)$. \\


\noindent Indeed,
since $u\in \mathcal{K}_\psi(\Omega)$, for every $v \ge 0$ and every $\e>0$ it results $u+\e v \ge \psi $, therefore if $v \in W^{1,q}_\mathrm{loc}(\Omega)$ by minimality of $u$

$$ \int_\Omega f(x, Du(x)) \,dx \le \int_\Omega f(x, Du+\varepsilon Dv(x))\, \dx$$
or equivalently

$$ \int_\Omega \Big[f(x, Du(x)+\varepsilon Dv(x))-f(x, Du(x))\Big]\,\dx \ge0.$$
Hence, we have

\begin{equation*}
 \varepsilon\int_\Omega \int_0^1\langle {D_{\xi}f}(x, Du(x)+\theta\varepsilon Dv(x)),Dv(x)\rangle\,\d\theta\,\dx \ge0
\end{equation*}
and also
\begin{equation*}
 \int_\Omega \int_0^1\langle {D_{\xi}f}(x, Du(x)+\theta\varepsilon Dv(x)),Dv(x)\rangle\,\d\theta\,\dx \ge0
\end{equation*}
where we divided both side of previous inequality by $\varepsilon$.
We observe that
\begin{eqnarray}\label{new}
0&\le& \int_\Omega \int_0^1\langle {D_{\xi}f}(x, Du(x)+\theta\varepsilon Dv),Dv\rangle\,\d\theta\,\dx \notag\\
&\le& \int_\Omega \int_0^1 | {D_{\xi}f}(x, Du+\theta\varepsilon Dv(x))||Dv(x)|\,\d\theta\,\dx \notag\\
&\le& \int_\Omega \int_0^1 (\mu^2+ |Du+\theta\varepsilon Dv(x)|^2)^{\frac{q-1}{2}}|Dv(x)|\,\d\theta\,\dx \notag\\
&\le& c\int_\Omega  (\mu^2+ |Du(x)|^2+\varepsilon^2 |Dv(x)|^2)^{\frac{q-1}{2}}|Dv(x)|\,\dx,
\end{eqnarray}
where in the last inequality we used Lemma 8.3 in \cite{Giusti}.\\
Therefore,  since $v \in W^{1,q}_\mathrm{loc}(\Omega)$, by the growth assumption (A3), assuming without loss of generality $\varepsilon<1$, we get

$$\int_0^1\langle {D_{\xi}f}(x, Du(x)+\theta\varepsilon Dv(x)),Dv(x)\rangle\,\d\theta \le \mu^q + |Du|^q+ |Dv|^q \in L^1(\Omega). $$
Then, applying dominated convergence theorem  in \eqref{new}, we have
$$\lim_{\varepsilon\to 0}\int_\Omega \int_0^1\langle {D_{\xi}f}(x, Du(x)+\theta\varepsilon Dv(x)),Dv(x)\rangle\,d\theta\,\dx =\int_\Omega \langle {D_{\xi}f}(x, Du(x)),Dv(x)\rangle\,\dx \ge0$$
for every $v\in W^{1,q}_0(\Omega)$, $v\ge 0$. At this point it is standard to verify
the inequality \eqref{variationalinequality}
$$\int_\Omega \langle {D_{\xi}f}(x, Du(x)),D\varphi(x)-Du(x)\rangle\,\dx\ge 0.$$

\noindent Now we have to choose suitable test functions $\varphi$ in
\eqref{variationalinequality} that involve the different quotient of
the solution and at the same time satisfy the conditions  $\varphi\in W^{1,q}_{\rm loc}(\Omega)$ and $\varphi\ge \psi$ in $\Omega$. In order to do this, we proceed similarly to what has been done
in \cite{CGG, EP}.

\bigskip

Let us fix a ball $B_{ R}\Subset \Omega$ and arbitrary radii $\frac{
R}{2}<r<s<t<\lambda r< R$, with $1<\lambda<2$. Let us consider a cut
off function $\eta\in C^\infty_0(B_t)$ such that $\eta\equiv 1$ on
$B_s$ and  $|D \eta|\le \frac{c}{t-s}$. From now on, with no
loss of generality, we suppose $R<1$.\\

Let $v\in W^{1,q}_0(\Omega)$ be such that

    \begin{equation}\label{cond}
    u-\psi+\tau v\ge 0\qquad\forall \tau\in[0, 1],
    \end{equation}
\noindent  and observe that $\varphi=u+\tau v\ge \psi $ for all $\tau \in [0, 1]$.
For $|h|<\frac{R}{4}$, we consider

    \begin{equation*}\label{v1}
    v_1(x)=\eta^2(x)\left[(u-\psi)(x+h)-(u-\psi)(x)\right],
    \end{equation*}

\noindent  so we have  $v_1\in W^{1, p+2}_0(\Omega)$, and,
for any $\tau\in[0,1]$, $v_1$ satisfies \eqref{cond}. Indeed, for
a.e. $x \in \Omega$ and for any $\tau\in[0,1]$

    \begin{eqnarray*}
    u(x)-\psi(x)+\tau v_1(x)&=&
    u(x)-\psi(x)+\tau\eta^2(x)\left[(u-\psi)(x+h)-(u-\psi)(x)\right] \cr\cr
    &=&\tau\eta^2(x)(u-\psi)(x+h)+(1-\tau\eta^2(x))(u-\psi)(x)\ge 0,
    \end{eqnarray*}

 \noindent   since $u\in \mathcal{K}_\psi(\Omega)$ and $0\le\eta\le1$. Therefore, from $q-p<1$ we have $L^{p+2}(\Omega)\hookrightarrow L^{q}(\Omega)$
    and so we can use $\varphi=u+\tau v_1$ as a test function in inequality \eqref{variationalinequality},  thus getting

    \begin{equation}\label{3.3}
    0\le\int_{\Omega}\left<A(x, Du(x)), D\left[\eta^2(x)\left[(u-\psi)(x+h)-(u-\psi)(x)\right]\right]\right>\dx.
    \end{equation}
Similarly, we define

    \begin{equation*}\label{v2}
    v_2(x)=\eta^2(x-h)\left[(u-\psi)(x-h)-(u-\psi)(x)\right],
    \end{equation*}
and we have $v_2\in W^{1, p+2}_0(\Omega)$, the inequality
\eqref{cond} still is satisfied for any $\tau\in[0,1]$, and we can use $\varphi=u+\tau v_2$ as test function in \eqref{variationalinequality},
obtaining


    \begin{equation*}
    0\le\int_{\Omega}\left<A(x, Du(x)), D\left[\eta^2(x-h)\left[(u-\psi)(x-h)-(u-\psi)(x)\right]\right]\right>\dx,
    \end{equation*}
and by means of a change of variable, we have

    \begin{equation}\label{3.5}
    0\le\int_{\Omega}\left<A(x+h, Du(x+h)), D\left[\eta^2(x)\left[(u-\psi)(x)-(u-\psi)(x+h)\right]\right]\right>\dx.
    \end{equation}

 \noindent  Now we can add \eqref{3.3} and \eqref{3.5}, thus getting

    \begin{align*}
    0\le&\int_{\Omega}\left<A(x, Du(x)), D\left[\eta^2(x)\left[(u-\psi)(x+h)-(u-\psi)(x)\right]\right]\right>\dx\\
    &+\int_{\Omega}\left<A(x+h, Du(x+h)), D\left[\eta^2(x)\left[(u-\psi)(x)-(u-\psi)(x+h)\right]\right]\right>\dx,
    \end{align*}
that is

\begin{equation*}
    0\le\int_{\Omega}\left<A(x, Du(x))-A(x+h, Du(x+h)), D\left[\eta^2(x)\left[(u-\psi)(x+h)-(u-\psi)(x)\right]\right]\right>\dx,
    \end{equation*}
which implies

    \begin{align*}
    0\ge&\int_{\Omega}\left<A(x+h, Du(x+h))-A(x, Du(x)), \eta^2(x)D\left[(u-\psi)(x+h)-(u-\psi)(x)\right]\right>\dx\notag\\
    &+\int_{\Omega}\left<A(x+h, Du(x+h))-A(x, Du(x)),
    2\eta(x)D\eta(x)\left[(u-\psi)(x+h)-(u-\psi)(x)\right]\right>\dx.
    \end{align*}

 \noindent  Previous inequality can be rewritten as follows

\begin{align}\label{differenceinequality3}
    0\ge&\int_{\Omega}\left<A(x+h, Du(x+h))-A(x+h, Du(x)),\eta^2(x)(Du(x+h)-Du(x))\right>\dx\notag\\
    &-\int_{\Omega}\left<A(x+h, Du(x+h))-A(x+h, Du(x)), \eta^2(x)(D\psi(x+h)-D\psi(x))\right>\dx\notag\\
    &+\int_{\Omega}\left<A(x+h, Du(x+h))-A(x+h, Du(x)), 2\eta(x)D\eta(x)\tau_h\left(u-\psi\right)(x)\right>\dx\notag\\
    &+\int_{\Omega}\left<A(x+h, Du(x))-A(x, Du(x)),\eta^2(x)(Du(x+h)-Du(x))\right>\dx\notag\\
    &-\int_{\Omega}\left<A(x+h, Du(x))-A(x, Du(x)), \eta^2(x)(D\psi(x+h)-D\psi(x))\right>\dx\notag\\
    &+\int_{\Omega}\left<A(x+h, Du(x))-A(x, Du(x)), 2\eta(x)D\eta(x)\tau_h\left(u-\psi\right)(x)\right>\dx\notag\\
    =:& \,I+II+III+IV+V+VI,
\end{align}

\noindent so we have

\begin{equation}\label{differenceinequality}
    I\le |II|+|III|+|IV|+|V|+|VI|.
\end{equation}

\noindent The ellipticity assumption (A1) implies

\begin{equation}\label{I}
    I\ge\nu\int_{\Omega}\eta^2(x)|\tau_hDu(x)|^2\left(\mu^2+\left|Du(x+h)\right|^2+\left|Du(x)\right|^2\right)^\frac{p-2}{2}\dx.
\end{equation}

\noindent By virtue of assumption (A2), Young's inequality with
exponents $\left(2, 2\right)$ and H\"{o}lder's inequality with
exponents $\left(\frac{p+2}{2(p+2-q)}, \frac{p+2}{2q-p-2}\right)$ we get

\begin{align*}
    |II|&\leq L\int_\om \eta^2(x)|\tau_hDu(x)|(\mu^2 +|Du(x)|^2+|Du(x+h)|^2)^\frac{q-2}{2}|\tau_hD\psi(x)|\dx\notag\\
	& \leq \e\int_\om \eta^2(x)|\tau_hDu(x)|^2(\mu^2 +|Du(x)|^2+|Du(x+h)|^2)^\frac{p-2}{2}\dx\notag\\
	& \quad \qquad + c_\e\int_\om \eta^2(x)|\tau_hD\psi(x)|^2(\mu^2 +|Du(x)|^2+|Du(x+h)|^2)^\frac{2q-p-2}{2}\dx\notag\\
	& \leq \e\int_\om \eta^2(x)|\tau_hDu(x)|^2(\mu^2 +|Du(x)|^2+|Du(x+h)|^2)^\frac{p-2}{2}\dx\notag\\
	& \quad \qquad + c_\e\left( \int_\bt  |\tau_hD\psi(x)|^\frac{p+2}{p+2-q} \dx\right)^\frac{2(p+2)-2q}{p+2} \cdot \left( \int_\bt(\mu^2 +|Du(x)|^2+|Du(x+h)|^2)^\frac{p+2}{2} \dx\right)^\frac{2q-p-2}{p+2}\notag\\
\end{align*}	
\noindent where we used also  the properties of $\eta$. Since $D\psi \in  W^{1,\frac{p+2}{p+2-q}}_\mathrm{loc}(\Omega)$, we may use the first and the second estimate of Lemma \ref{le1}
to control the first and the second integral respectively in the last line of the previous estimate, thus obtaining

\begin{align}\label{II}
 |II| & \leq \e\int_\om \eta^2(x)|\tau_hDu(x)|^2(\mu^2 +|Du(x)|^2+|Du(x+h)|^2)^\frac{p-2}{2}\dx\notag\\
	  & \quad \qquad+ c_\e|h|^2\left(\int_\blr |D^2\psi(x)|^\frac{p+2}{p+2-q} \dx\right)^\frac{2(p+2-q)}{p+2} \cdot \left(\int_\blr (\mu^2 +|Du(x)|^2)^\frac{p+2}{2} \dx\right)^\frac{2q-p-2}{p+2}.
\end{align}

\medskip

\noindent Arguing analogously, by virtue of assumption (A2) and Young's and H\"{o}lder's inequalities, we get

\begin{align*}
|III|&\le 2L\int_{\Omega}\eta(x)|D \eta(x)||\tau_hDu(x)|\left(\mu^2+|Du(x+h)|^2+|Du(x)|^2\right)^\frac{q-2}{2}|\tau_h\left(u-\psi\right)|\dx\notag\\
     & \le \e\int_\om \eta^2(x)|\tau_hDu(x)|^2(\mu^2 +|Du(x)|^2+|Du(x+h)|^2)^\frac{p-2}{2}\dx\\
	 &\quad \qquad + c_\e \int_\om |\tau_h(u-\psi)(x)|^2|D\eta(x)|^2(\mu^2 +|Du(x)|^2+|Du(x+h)|^2)^\frac{2q-p-2}{2}\dx\\
	& \leq \e\int_\om \eta^2(x)|\tau_hDu(x)|^2(\mu^2 +|Du(x)|^2+|Du(x+h)|^2)^\frac{p-2}{2}\dx\\
	& \quad \qquad + \frac{c_\e}{(t-s)^2} \int_{B_t} |\tau_h(u-\psi)(x)|^2(\mu^2 +|Du(x)|^2+|Du(x+h)|^2)^\frac{2q-p-2}{2}\dx\\
	& \leq \e\int_\om \eta^2(x)|\tau_hDu(x)|^2(\mu^2 +|Du(x)|^2+|Du(x+h)|^2)^\frac{p-2}{2}\dx\notag\\
	& \quad \qquad + \frac{c_\e}{(t-s)^2} \left( \int_\bt |\tau_h(u-\psi)(x)|^\frac{p+2}{p+2-q} \dx\right)^\frac{2(p+2)-2q}{p+2} \cdot \\
    & \qquad \qquad \qquad \qquad \cdot \left( \int_\bt(\mu^2 +|Du(x)|^2+|Du(x+h)|^2)^\frac{p+2}{2} \dx\right)^\frac{2q-p-2}{p+2}\notag
\end{align*}
\noindent where we used the properties of $\eta$. Using now both estimates of Lemma \ref{le1}, we get

\begin{align}\label{III}
|III|& \leq \e\int_\om \eta^2(x)|\tau_hDu(x)|^2(\mu^2 +|Du(x)|^2+|Du(x+h)|^2)^\frac{p-2}{2}\dx \notag\\
	& \quad \qquad + \frac{c_\e|h|^2}{(t-s)^2} \left(\int_\blr |D(u-\psi)|^\frac{p+2}{p+2-q}\dx\right)^\frac{2(p+2)-2q}{p+2} \cdot \left( \int_\blr \mu^{p+2} +|Du(x)|^{p+2}\dx \right)^\frac{2q-p-2}{p+2}.
\end{align}

\noindent In order to estimate the term $IV$, 
we use assumption (A4) and Young's inequality to obtain

\begin{align*}
|IV|& \le |h|\int_{\Omega}\eta^2(x)\left(\kappa(x+h)+\kappa(x)\right)\left(\mu^2+|Du(x)|^2\right)^\frac{q-1}{2}|\tau_hDu(x)|\dx\notag\\
& \le \varepsilon\int_{\Omega}\eta^2(x)\left|\tau_hDu(x)\right|^2\left(\mu^2+|Du(x+h)|^2+|Du(x)|^2\right)^\frac{p-2}{2}\dx\notag\\
&\quad \qquad +c_\varepsilon|h|^2\int_{B_t}\left(\kappa(x+h)+\kappa(x)\right)^2\left(\mu^2+|Du(x)|^2\right)^\frac{2q-p}{2}\dx.
\end{align*}
Using H\"{o}lder's inequality with exponents
$\left(\frac{p+2}{2(p-q+1)}, \frac{p+2}{2q-p}\right)$ and Lemma \ref{le1} we have

\begin{align}\label{IV}
|IV|& \le\varepsilon\int_{\Omega}\eta^2(x)\left|\tau_hDu(x)\right|^2\left(\mu^2+|Du(x+h)|^2+|Du(x)|^2\right)^\frac{p-2}{2}\dx\notag\\
&\quad \qquad +c_\varepsilon|h|^2\left(\int_{B_{\lambda r}}\kappa^{\frac{p+2}{p-q+1}}(x)\dx\right)^\frac{2p-2q+2}{p+2}\cdot\left(\int_{B_t}\left(\mu^{p+2}+|Du(x)|^{p+2}\right)\dx\right)^\frac{2q-p}{p+2}.
\end{align}

The condition (A4) also entails

\begin{align}\label{V}
|V|& \le|h|\int_{\Omega}\eta^2(x)\left(\kappa(x+h)+\kappa(x)\right)\left(\mu^2+|Du(x)|^2\right)^\frac{q-1}{2}\left|\tau_hD\psi(x)\right|\dx\notag\\
&\le |h|\left(\int_{B_t}\left(\kappa(x+h)+\kappa(x)\right)^{\frac{p+2}{p-q+1}}\dx\right)^\frac{p-q+1}{p+2}\cdot\left(\int_{B_t}\left(\mu^{p+2}+|Du(x)|^{p+2}\right)dx\right)^\frac{q-1}{p+2}\notag\\
&\quad \qquad \cdot\left(\int_{B_t}\left|\tau_hD\psi(x)\right|^\frac{p+2}{2}\dx\right)^\frac{2}{p+2}\notag\\
&\le c|h|^2\left(\int_{B_{\lambda r}}\kappa^{\frac{p+2}{p-q+1}}(x)\dx\right)^\frac{p-q+1}{p+2}\cdot\left(\int_{B_t}\left(\mu^{p+2}+|Du(x)|^{p+2}\right)\dx\right)^\frac{q-1}{p+2}\cdot\left(\int_{B_{\lambda r}}\left|D^2\psi(x)\right|^\frac{p+2}{2}\dx\right)^\frac{2}{p+2}
\end{align}
where we used H\"{o}lder's inequality with exponents $\left(\frac{p+2}{p-q+1},
\frac{p+2}{q-1}, \frac{p+2}{2}\right)$, the properties of $\eta$ and  Lemma \ref{le1}.

 Finally, using again assumption (A4), the properties of $\eta$, H\"{o}lder's inequality and Lemma \ref{le1}, we have

\begin{align}\label{VI}
|VI|& \le 2|h|\int_{\Omega}\eta(x)\left|D\eta(x)\right|\left(\kappa(x+h)+\kappa(x)\right)\left(\mu^2+|Du(x)|^2\right)^\frac{q-1}{2}
\left|\tau_h\left(u-\psi\right)(x)\right|\dx\notag\\
& \le \frac{c|h|}{t-s}\left(\int_{B_t}\left(\kappa(x+h)+\kappa(x)\right)^{\frac{p+2}{p-q+1}}\dx\right)^\frac{p-q+1}{p+2}\cdot\left(\int_{B_t}\left(\mu^{p+2}+|Du(x)|^{p+2}\right)\dx\right)^\frac{q-1}{p+2}\notag\\
&\quad \qquad \cdot\left(\int_{B_t}\left|\tau_h\left(u-\psi\right)(x)\right|^\frac{p+2}{2}\dx\right)^\frac{2}{p+2}\notag\\
\le&\frac{c|h|^2}{t-s}\left(\int_{B_{\lambda r}}\kappa(x)^{\frac{p+2}{p-q+1}}\dx\right)^\frac{p-q+1}{p+2}\cdot\left(\int_{B_t}\left(\mu^{p+2}+|Du(x)|^{p+2}\right)dx\right)^\frac{q-1}{p+2}\notag\\
&\quad \qquad\cdot\left(\int_{B_{\lambda r}}\left|D\left(u-\psi\right)(x)\right|^\frac{p+2}{2}\dx\right)^\frac{2}{p+2}.
\end{align}

Inserting \eqref{I}, \eqref{II}, \eqref{III}, \eqref{IV}, \eqref{V}
and \eqref{VI} into \eqref{differenceinequality} we infer
\begin{align*}
\nu\int_{\Omega}&\eta^2(x)|\tau_hDu(x)|^2\left(\mu^2+\left|Du(x+h)\right|^2+\left|Du(x)\right|^2\right)^\frac{p-2}{2}\dx\notag\\
& \leq 3\e\int_\om \eta(x)^2|\tau_hDu(x)|^2(\mu^2 +|Du(x)|^2+|Du(x+h)|^2)^\frac{p-2}{2}\dx\notag\\
&	\qquad+ c_\e|h|^2\left(\int_\blr |D^2\psi(x)|^\frac{p+2}{p+2-q} \dx\right)^\frac{2(p+2)-2q}{p+2} \cdot \left(\int_\blr (\mu^2 +|Du(x)|^2)^\frac{p+2}{2} \dx\right)^\frac{2q-p-2}{p+2}\notag\\
&  \qquad + \frac{c_\e|h|^2}{(t-s)^2} \left(\int_\blr |D(u-\psi)|^\frac{p+2}{p+2-q}\dx\right)^\frac{2(p+2)-2q}{p+2} \cdot \left( \int_\blr \mu^{p+2} +|Du(x)|^{p+2}\dx \right)^\frac{2q-p-2}{p+2}\notag\\
& \qquad +c_\varepsilon|h|^2\left(\int_{B_{\lambda r}}\kappa^{\frac{p+2}{p-q+1}}(x)\dx\right)^\frac{2p-2q+2}{p+2}\cdot\left(\int_{B_t}\left(\mu^{p+2}+|Du(x)|^{p+2}\right)\dx\right)^\frac{2q-p}{p+2}\notag\\
&\qquad +c|h|^2\left(\int_{B_{\lambda r}}\kappa^{\frac{p+2}{p-q+1}}(x)\dx\right)^\frac{p-q+1}{p+2}\cdot\left(\int_{B_t}\left(\mu^{p+2}+|Du(x)|^{p+2}\right)\dx\right)^\frac{q-1}{p+2}\notag\\
&\quad \qquad \cdot\left(\int_{B_{\lambda r}}\left|D^2\psi(x)\right|^\frac{p+2}{2}\dx\right)^\frac{2}{p+2}\notag\\
&\qquad +\frac{c|h|^2}{t-s}\left(\int_{B_{\lambda r}}\kappa(x)^{\frac{p+2}{p-q+1}}\dx\right)^\frac{p-q+1}{p+2}\cdot\left(\int_{B_t}\left(\mu^{p+2}+|Du(x)|^{p+2}\right)\dx\right)^\frac{q-1}{p+2}\notag\\
&\quad \qquad\cdot\left(\int_{B_{\lambda r}}\left|D\left(u-\psi\right)(x)\right|^\frac{p+2}{2}\dx\right)^\frac{2}{p+2}.
\end{align*}


Choosing
$\varepsilon=\frac{\nu}{6}$, we can reabsorb the first term from the
right-hand side to the left-hand one, thus getting

\begin{align}\label{Plugestimates}
\nu &\int_{\Omega}\eta^2(x)|\tau_hDu(x)|^2\left(\mu^2+\left|Du(x+h)\right|^2+\left|Du(x)\right|^2\right)^\frac{p-2}{2}\dx\notag\\
&\le c|h|^2\left(\int_\blr |D^2\psi(x)|^\frac{p+2}{p+2-q} \dx\right)^\frac{2(p+2)-2q}{p+2} \cdot \left(\int_\blr (\mu^2 +|Du(x)|^2)^\frac{p+2}{2} \dx\right)^\frac{2q-p-2}{p+2}\notag\\
&  \quad + \frac{c|h|^2}{(t-s)^2} \left(\int_\blr |D(u-\psi)|^\frac{p+2}{p+2-q}\dx\right)^\frac{2(p+2)-2q}{p+2} \cdot \left( \int_\blr \mu^{p+2} +|Du(x)|^{p+2}\dx \right)^\frac{2q-p-2}{p+2}\notag\\
& \quad +c|h|^2\left(\int_{B_{\lambda r}}\kappa^{\frac{p+2}{p-q+1}}(x)\dx\right)^\frac{2p-2q+2}{p+2}\cdot\left(\int_{B_t}\left(\mu^{p+2}+|Du(x)|^{p+2}\right)\dx\right)^\frac{2q-p}{p+2}\notag\\
&\quad +c |h|^2\left(\int_{B_{\lambda r}}\kappa^{\frac{p+2}{p-q+1}}(x)\dx\right)^\frac{p-q+1}{p+2}\cdot\left(\int_{B_t}\left(\mu^{p+2}+|Du(x)|^{p+2}\right)\dx\right)^\frac{q-1}{p+2}\cdot\left(\int_{B_{\lambda r}}\left|D^2\psi(x)\right|^\frac{p+2}{2}\dx\right)^\frac{2}{p+2}\notag\\
&\quad +\frac{c|h|^2}{t-s}\left(\int_{B_{\lambda r}}\kappa(x)^{\frac{p+2}{p-q+1}}\dx\right)^\frac{p-q+1}{p+2}\cdot\left(\int_{B_t}\left(\mu^{p+2}+|Du(x)|^{p+2}\right)\dx\right)^\frac{q-1}{p+2}\notag\\
&\quad \qquad\cdot\left(\int_{B_{\lambda r}}\left|D\left(u-\psi\right)(x)\right|^\frac{p+2}{2}\dx\right)^\frac{2}{p+2}.
\end{align}


\noindent Now we apply Young's inequalities  and since $u\in \mathcal{K}_{\psi}(\Omega)$,
we have

\begin{align}\label{Plugestimates}
\nu & \int_{\Omega}\eta^2(x)|\tau_hDu(x)|^2\left(\mu^2+\left|Du(x+h)\right|^2+\left|Du(x)\right|^2\right)^\frac{p-2}{2}\dx\notag\\
&	\le  c_\e|h|^2 \int_\blr |D^2\psi(x)|^\frac{p+2}{p+2-q} \dx + \e|h|^2 \int_\blr (\mu^{p+2} +|Du(x)|^{p+2}) \dx\notag\\
&  \quad + \frac{c_\e|h|^2}{(t-s)^{\frac{p+2}{p-q+2}}} \int_\blr |D(u-\psi)|^\frac{p+2}{p+2-q}\dx + \e |h|^2  \int_\blr (\mu^{p+2} +|Du(x)|^{p+2})\dx \notag\\
& \quad +c_\e|h|^2\int_{B_{\lambda r}}\kappa^{\frac{p+2}{p-q+1}}(x)\dx + \e |h|^2  \int_{B_t}\left(\mu^{p+2}+|Du(x)|^{p+2}\right)\dx\notag\\
&\quad +c_\e |h|^2 \int_{B_{\lambda r}}\kappa^{\frac{p+2}{p-q+1}}(x)\dx + \e |h|^2 \int_{B_t}\left(\mu^{p+2}+|Du(x)|^{p+2}\right)\dx + c_\e |h|^2 \int_{B_{\lambda r}}\left|D^2\psi(x)\right|^\frac{p+2}{2}\dx \notag\\
&\quad  + c_\e |h|^2 \int_{B_{\lambda r}}\kappa(x)^{\frac{p+2}{p-q+1}}\dx+ \e |h|^2 \int_{B_t}\left(\mu^{p+2}+|Du(x)|^{p+2}\right)\dx\notag\\
&\quad +\frac{c_\e|h|^2}{(t-s)^{\frac{p+2}{2}}} \int_{B_{\lambda r}}\left|D\left(u-\psi\right)(x)\right|^\frac{p+2}{2}\dx.
\end{align}

By Young's inequalities of exponents $(p+2-q, \frac{p+2-q}{p+1-q})$ and $(2,2)$ we can estimate the third and the last integral appearing in the right hand side  of the previous inequality as

\begin{align*}
	& \frac{c_\e|h|^2}{(t-s)^{\frac{p+2}{p-q+2}}} \int_\blr |D(u-\psi)|^\frac{p+2}{p+2-q}\,\dx\\
	&\qquad \leq \frac{c_\e|h|^2}{(t-s)^{\frac{p+2}{p-q+2}}} \int_\blr |Du(x)|^\frac{p+2}{p+2-q}\dx+\frac{c_\e|h|^2}{(t-s)^{\frac{p+2}{p-q+2}}}\int_\blr|D\psi(x)|^\frac{p+2}{p+2-q}\, \dx\\
    & \qquad\le c_\e|h|^2 \int_\blr  |Du(x)|^{p+2} \,\dx + \frac{c_\e|h|^2}{(t-s)^{\frac{p+2}{p-q+1}}}|B_R| +\frac{c_\e|h|^2}{(t-s)^{\frac{p+2}{p-q+2}}}\int_\blr|D\psi(x)|^\frac{p+2}{p+2-q}\, \dx,\\
	& \qquad\leq \frac{c_\e|h|^2}{(t-s)^{\frac{p+2}{p-q+1}}}|B_R| +\e|h|^2 \int_\blr (\mu^{p+2}+ |Du(x)|^{p+2}) \, \dx+\frac{c_\e|h|^2}{(t-s)^{\frac{p+2}{p-q+2}}}\int_\blr|D\psi(x)|^\frac{p+2}{p+2-q}\, \dx,
\end{align*}

and similarly

\begin{align*}
	& \frac{c_\e|h|^2}{(t-s)^{\frac{p+2}{2}}} \int_\blr |D(u-\psi)|^\frac{p+2}{2}\,\dx\\
		& \qquad\leq \frac{c_\e|h|^2}{(t-s)^{p+2}}|B_R| +\e|h|^2 \int_\blr (\mu^{p+2}+ |Du(x)|^{p+2}) \, \dx+\frac{c_\e|h|^2}{(t-s)^{\frac{p+2}{2}}}\int_\blr|D\psi(x)|^\frac{p+2}{2}\, \dx.
\end{align*}

So, from \eqref{Plugestimates}, we get

\begin{align}\label{Plugestimates2}
\nu\int_{\Omega}&\eta^2(x)|\tau_hDu(x)|^2\left(\mu^2+\left|Du(x+h)\right|^2+\left|Du(x)\right|^2\right)^\frac{p-2}{2}\dx\notag\\
&	\le  c_\e|h|^2 \int_\blr |D^2\psi(x)|^\frac{p+2}{p+2-q} \dx + \e|h|^2 \int_\blr (\mu^{p+2} +|Du(x)|^{p+2})\,  \dx\notag\\
&  \qquad +\frac{c_\e|h|^2}{(t-s)^{\frac{p+2}{p-q+1}}}|B_R| +2\e|h|^2 \int_\blr (\mu^{p+2}+ |Du(x)|^{p+2}) \, \dx+\frac{c_\e|h|^2}{(t-s)^{\frac{p+2}{p-q+2}}}\int_\blr|D\psi(x)|^\frac{p+2}{p+2-q}\, \dx, \notag\\
& \qquad +c_\e|h|^2\int_{B_{\lambda r}}\kappa^{\frac{p+2}{p-q+1}}(x)\dx + \e |h|^2  \int_{B_t}\left(\mu^{p+2}+|Du(x)|^{p+2}\right)\dx\notag\\
&\qquad +c_\e |h|^2 \int_{B_{\lambda r}}\kappa^{\frac{p+2}{p-q+1}}(x)\dx + \e |h|^2 \int_{B_t}\left(\mu^{p+2}+|Du(x)|^{p+2}\right)\dx + c_\e |h|^2 \int_{B_{\lambda r}}\left|D^2\psi(x)\right|^\frac{p+2}{2}\, \dx \notag\\
&\qquad  + c_\e |h|^2 \int_{B_{\lambda r}}\kappa(x)^{\frac{p+2}{p-q+1}}\dx+ \e |h|^2 \int_{B_t}\left(\mu^{p+2}+|Du(x)|^{p+2}\right)\dx\notag\\
&\quad \qquad +\frac{c_\e|h|^2}{(t-s)^{p+2}}|B_R| +\e|h|^2 \int_\blr (\mu^{p+2}+ |Du(x)|^{p+2}) \, \dx+\frac{c_\e|h|^2}{(t-s)^{\frac{p+2}{2}}}\int_\blr|D\psi(x)|^\frac{p+2}{2}\, \dx\notag\\
&	\le  c_\e|h|^2 \int_\blr |D^2\psi(x)|^\frac{p+2}{p+2-q} \dx + C \cdot\e|h|^2 \int_\blr (\mu^{p+2} +|Du(x)|^{p+2}) \, \dx\notag\\
& \qquad +\frac{c_\e|h|^2}{(t-s)^{\frac{p+2}{p-q+1}}}|B_R| +\frac{c_\e|h|^2}{(t-s)^{\frac{p+2}{p-q+2}}}\int_\blr|D\psi(x)|^\frac{p+2}{p+2-q}\, \dx +3c_\e |h|^2 \int_{B_{\lambda r}}\kappa^{\frac{p+2}{p-q+1}}(x)\dx \notag\\
&\qquad   + c_\e |h|^2 \int_{B_{\lambda r}}\left|D^2\psi(x)\right|^\frac{p+2}{2}\, \dx +\frac{c_\e|h|^2}{(t-s)^{p+2}}|B_R| +\frac{c_\e|h|^2}{(t-s)^{\frac{p+2}{2}}}\int_\blr|D\psi(x)|^\frac{p+2}{2}\, \dx
\end{align}

\noindent Using, in the left hand side of the previous estimate, the right-hand side of the inequality
\eqref{lemma6GPestimate1} in Lemma \ref{lemma6GP} , we get

\begin{align}
\nu\int_\Omega&\eta^2(x)\left|\tau_hV_p\left(Du(x)\right)\right|^2\dx\notag\\
&	\le  c_\e|h|^2 \int_\blr |D^2\psi(x)|^\frac{p+2}{p+2-q} \dx + 7\e|h|^2 \int_\blr (\mu^{p+2} +|Du(x)|^{p+2}) \, \dx\notag\\
& \qquad +\frac{c_\e|h|^2}{(t-s)^{\frac{p+2}{p-q+1}}}|B_R| +\frac{c_\e|h|^2}{(t-s)^{\frac{p+2}{p-q+2}}}\int_\blr|D\psi(x)|^\frac{p+2}{p+2-q}\, \dx +3c_\e |h|^2 \int_{B_{\lambda r}}\kappa^{\frac{p+2}{p-q+1}}(x)\dx \notag\\
&\qquad   + c_\e |h|^2 \int_{B_{\lambda r}}\left|D^2\psi(x)\right|^\frac{p+2}{2}\, \dx +\frac{c_\e|h|^2}{(t-s)^{p+2}}|B_R| +\frac{c_\e|h|^2}{(t-s)^{\frac{p+2}{2}}}\int_\blr|D\psi(x)|^\frac{p+2}{2}\, \dx
\end{align}

\noindent Dividing both sides by $|h|^2$ and using Lemma
\ref{Giusti8.2}, we have

\begin{align}\label{stima derivate seconde}
    \int_{B_t}&\left|DV_p(Du(x))\right|^2\dx\le  c_\e\int_\blr |D^2\psi(x)|^\frac{p+2}{p+2-q} \dx + 7\e \int_\blr (\mu^{p+2} +|Du(x)|^{p+2}) \, \dx\notag\\
& \qquad +\frac{c_\e}{(t-s)^{\frac{p+2}{p-q+1}}}|B_R| +\frac{c_\e}{(t-s)^{\frac{p+2}{p-q+2}}}\int_\blr|D\psi(x)|^\frac{p+2}{p+2-q}\, \dx +3c_\e  \int_{B_{\lambda r}}\kappa^{\frac{p+2}{p-q+1}}(x)\dx \notag\\
&\qquad   + c_\e  \int_{B_{\lambda r}}\left|D^2\psi(x)\right|^\frac{p+2}{2}\, \dx +\frac{c_\e}{(t-s)^{p+2}}|B_R| +\frac{c_\e}{(t-s)^{\frac{p+2}{2}}}\int_\blr|D\psi(x)|^\frac{p+2}{2}\, \dx.
\end{align}

\noindent Now, by virtue of left-hand side of inequality
\eqref{lemma6GPestimate2} of Lemma \ref{lemma6GP}

\begin{align}\label{beforeInterpolation}
\int_{B_t}&\left(\mu^2+\left|Du(x)\right|^2\right)^\frac{p-2}{2}\left|D^2u(x)\right|^2\dx\le\int_{B_t}\left|DV_p(Du(x))\right|^2\dx\notag\\
&\le c_\e\int_\blr |D^2\psi(x)|^\frac{p+2}{p+2-q} \dx + 7\e \int_\blr (\mu^{p+2} +|Du(x)|^{p+2}) \, \dx\notag\\
& \qquad +\frac{c_\e}{(t-s)^{\frac{p+2}{p-q+1}}}|B_R| +\frac{c_\e}{(t-s)^{\frac{p+2}{p-q+2}}}\int_\blr|D\psi(x)|^\frac{p+2}{p+2-q}\, \dx +3c_\e  \int_{B_{\lambda r}}\kappa^{\frac{p+2}{p-q+1}}(x)\dx \notag\\
&\qquad   + c_\e  \int_{B_{\lambda r}}\left|D^2\psi(x)\right|^\frac{p+2}{2}\, \dx +\frac{c_\e}{(t-s)^{p+2}}|B_R| +\frac{c_\e}{(t-s)^{\frac{p+2}{2}}}\int_\blr|D\psi(x)|^\frac{p+2}{2}\, \dx.
\end{align}

By virtue of the local boundedness of $u$, the second interpolation inequality of Lemma \ref{lemma5GP} yields that

\begin{align*}
\int_{\Omega}&\eta^2(x)\left(\mu^2+\left|Du(x)\right|^2\right)^\frac{p}{2}\left|Du(x)\right|^2\dx\notag\\
\le&c\Arrowvert u\Arrowvert_{L^\infty\left(\mathrm{supp}(\eta)\right)}^2\int_\Omega\eta^2(x)\left(\mu^2+\left|Du(x)\right|^2\right)^\frac{p-2}{2}\left|D^2u(x)\right|^2\dx\notag\\
&+c\Arrowvert
u\Arrowvert_{L^\infty\left(\mathrm{supp}(\eta)\right)}^2\int_{\Omega}\left(|\eta(x)|^2+\left|D\eta(x)\right|^2\right)
\left(\mu^2+\left|Du(x)\right|^2\right)^\frac{p}{2}\dx.
\end{align*}

\noindent and so,  combining  this last  estimate with \eqref{beforeInterpolation}, and using the properties of $\eta$, we get

\begin{align}\label{beforeIteration1}
&\int_{B_{r}}\left(\mu^2+\left|Du(x)\right|^2\right)^\frac{p}{2}\left|Du(x)\right|^2 \, \dx \notag\\
\le & c_\e \Arrowvert u\Arrowvert_{L^\infty\left(B_{\lambda r}\right)}^2\int_\blr |D^2\psi(x)|^\frac{p+2}{p+2-q} \dx + 7\e \Arrowvert u\Arrowvert_{L^\infty\left(B_{\lambda r}\right)}^2\int_\blr (\mu^{p+2} +|Du(x)|^{p+2}) \, \dx\notag\\
& +\frac{c_\e \Arrowvert u\Arrowvert_{L^\infty\left(B_{\lambda r}\right)}^2}{(t-s)^{\frac{p+2}{p-q+1}}}|B_R| +\frac{c_\e \Arrowvert u\Arrowvert_{L^\infty\left(B_{\lambda r}\right)}^2}{(t-s)^{\frac{p+2}{p-q+2}}}\int_\blr|D\psi(x)|^\frac{p+2}{p+2-q}\, \dx +3c_\e \Arrowvert u\Arrowvert_{L^\infty\left(B_{\lambda r}\right)}^2  \int_{B_{\lambda r}}\kappa^{\frac{p+2}{p-q+1}}(x)\dx \notag\\
&  + c_\e \Arrowvert u\Arrowvert_{L^\infty\left(B_{\lambda r}\right)}^2 \int_{B_{\lambda r}}\left|D^2\psi(x)\right|^\frac{p+2}{2}\, \dx +\frac{c_\e \Arrowvert u\Arrowvert_{L^\infty\left(B_{\lambda r}\right)}^2}{(t-s)^{p+2}}|B_R| +\frac{c_\e \Arrowvert u\Arrowvert_{L^\infty\left(B_{\lambda r}\right)}^2}{(t-s)^{\frac{p+2}{2}}}\int_\blr|D\psi(x)|^\frac{p+2}{2}\, \dx\notag\\
&+\frac{c_\varepsilon\Arrowvert u\Arrowvert_{L^\infty\left(B_{\lambda r}\right)}^2}{(t-s)^2}\int_{B_{\lambda r}}\left(\mu^2+\left|Du(x)\right|^2\right)^\frac{p}{2}\, \dx,
\end{align}
that we can also rewrite as

\begin{align*}
&\int_{B_{r}} \left(\mu^2+\left|Du(x)\right|^2\right)^\frac{p}{2}\left|Du(x)\right|^2 \dx\notag\\
\le & 7\e \Arrowvert u\Arrowvert_{L^\infty\left(B_{\lambda r}\right)}^2\int_\blr (\mu^{p+2} +|Du(x)|^{p+2}) \, \dx + c_\e \Arrowvert u\Arrowvert_{L^\infty\left(B_{\lambda r}\right)}^2\int_\blr |D^2\psi(x)|^\frac{p+2}{p+2-q} \dx \notag\\
& +c_\e \Arrowvert u\Arrowvert_{L^\infty\left(B_{\lambda r}\right)}^2  \int_{B_{\lambda r}}\kappa^{\frac{p+2}{p-q+1}}(x)\dx  + c_\e \Arrowvert u\Arrowvert_{L^\infty\left(B_{\lambda r}\right)}^2 \int_{B_{\lambda r}}\left|D^2\psi(x)\right|^\frac{p+2}{2}\, \dx \notag\\
& +\frac{c_\e \Arrowvert u\Arrowvert_{L^\infty\left(B_{\lambda r}\right)}^2}{(t-s)^{\frac{p+2}{p-q+1}}}\left[|B_R| +\int_\blr|D\psi(x)|^\frac{p+2}{p+2-q}\, \dx +\int_\blr|D\psi(x)|^\frac{p+2}{2}\, \dx
+\int_{B_{\lambda r}}\left(\mu^2+\left|Du(x)\right|^2\right)^\frac{p}{2}\, \dx \right]
\end{align*}
since $p\ge2$ and  $t-s<1$. 

\noindent Choosing $\varepsilon$ such that $7 \varepsilon\ \Arrowvert
u\Arrowvert_{L^\infty\left(B_{R}\right)}^2\le\frac{1}{2}$, previous
estimate becomes

\begin{align}\label{beforeIteration3}
\int_{B_{r}}&\left|Du(x)\right|^{p+2} \dx\le\int_{B_{r}}\left(\mu^2+\left|Du(x)\right|^2\right)^\frac{p}{2}\left|Du(x)\right|^2 \, \dx\le
\frac{1}{2}\int_{B_{\lambda r}}\left|Du(x)\right|^{p+2}\, \dx +c( \mu,p)|B_{R}|\notag\\
&+c\Arrowvert
u\Arrowvert_{L^\infty\left(B_{R}\right)}^2\left[   \int_{B_{R}} |D^2\psi(x)|^\frac{p+2}{p+2-q} \dx  +\int_{B_{R}}\kappa^{\frac{p+2}{p-q+1}}(x)\dx +    \int_{B_{R}}\left|D^2\psi(x)\right|^\frac{p+2}{2}\, \dx \right]\notag\\
& +\frac{c \Arrowvert u\Arrowvert_{L^\infty\left(B_{R}\right)}^2}{(t-s)^{\frac{p+2}{p-q+1}}}\left[|B_R| +\int_{B_{R}}|D\psi(x)|^\frac{p+2}{p+2-q}\, \dx +\int_{B_{R}}|D\psi(x)|^\frac{p+2}{2}\, \dx
+\int_{B_{R}}\left(\mu^2+\left|Du(x)\right|^2\right)^\frac{p}{2}\, \dx \right]
\end{align}

%

\noindent where $c=c(p, q, L, \nu, \mu)$ is independent of $t$ and $s$. Since
\eqref{beforeIteration3} is valid for any $\frac{R}{2}<r<s<t<\lambda
r<R<1$, taking the limit as $s\to r$ and $t\to\lambda r$, we get

\begin{align}\label{beforeIteration5}
&\int_{B_{r}}\left|Du(x)\right|^{p+2} \dx\le
\frac{1}{2}\int_{B_{\lambda r}}\left|Du(x)\right|^{p+2}\, \dx +c( \mu,p)|B_{R}|\notag\\
+&c\Arrowvert
u\Arrowvert_{L^\infty\left(B_{R}\right)}^2\left[   \int_{B_{R}} |D^2\psi(x)|^\frac{p+2}{p+2-q} \dx  +\int_{B_{R}}\kappa^{\frac{p+2}{p-q+1}}(x)\dx +    \int_{B_{R}}\left|D^2\psi(x)\right|^\frac{p+2}{2}\, \dx \right]\notag\\
+ &\frac{c \Arrowvert u\Arrowvert_{L^\infty\left(B_{R}\right)}^2}{r^{\frac{p+2}{p-q+1}} (\lambda-1)^{\frac{p+2}{p-q+1}}}\left[|B_R| +\int_{B_{R}}|D\psi(x)|^\frac{p+2}{p+2-q}\, \dx +\int_{B_{R}}|D\psi(x)|^\frac{p+2}{2}\, \dx
+\int_{B_{R}}\left(\mu^2+\left|Du(x)\right|^2\right)^\frac{p}{2}\, \dx \right]
\end{align}

Now, setting

\begin{equation*}
 h(r)=\int_{B_{r}}\left|Du(x)\right|^{p+2} \dx,
\end{equation*}

\begin{equation*}
 A=c\Arrowvert u\Arrowvert_{L^\infty\left(B_{R}\right)}^2\left[|B_R| +\int_{B_{R}}|D\psi(x)|^\frac{p+2}{p+2-q}\, \dx +\int_{B_{R}}|D\psi(x)|^\frac{p+2}{2}\, \dx
+\int_{B_{R}}\left(\mu^2+\left|Du(x)\right|^2\right)^\frac{p}{2}\, \dx \right],
\end{equation*}

and

\begin{equation*}
B=c( \mu,p)|B_{R}|+c\Arrowvert
u\Arrowvert_{L^\infty\left(B_{R}\right)}^2\left[   \int_{B_{R}} |D^2\psi(x)|^\frac{p+2}{p+2-q} \dx  +\int_{B_{R}}\kappa^{\frac{p+2}{p-q+1}}(x)\dx +    \int_{B_{R}}\left|D^2\psi(x)\right|^\frac{p+2}{2}\, \dx\right],
\end{equation*}
we obtain

\begin{equation*}
 h(r)\le \frac 12 h(\lambda r) + \frac{A}{r^{\frac{p+2}{p-q+1}} (\lambda-1)^{\frac{p+2}{p-q+1}}}+B
\end{equation*}

Thus, we can apply Lemma \ref{iteration}, with

\begin{equation*}
\theta=\frac{1}{2}\qquad\text{ and }\qquad\gamma=\frac{p+2}{p-q+1},
\end{equation*}

obtaining

\begin{align*}
&\int_{B_{r}}\left|Du(x)\right|^{p+2}\, \dx  \notag\\
\le &c( \mu,p)|B_{R}|+c\Arrowvert
u\Arrowvert_{L^\infty\left(B_{R}\right)}^2\left[   \int_{B_{R}} |D^2\psi(x)|^\frac{p+2}{p+2-q} \dx  +\int_{B_{R}}\kappa^{\frac{p+2}{p-q+1}}(x)\dx +    \int_{B_{R}}\left|D^2\psi(x)\right|^\frac{p+2}{2}\, \dx \right]\notag\\
& +\frac{c \Arrowvert u\Arrowvert_{L^\infty\left(B_{R}\right)}^2}{ R^{\frac{p+2}{p-q+1}}}\left[|B_R| +\int_{B_{R}}|D\psi(x)|^\frac{p+2}{p+2-q}\, \dx +\int_{B_{R}}|D\psi(x)|^\frac{p+2}{2}\, \dx
+\int_{B_{R}}\left(\mu^2+\left|Du(x)\right|^2\right)^\frac{p}{2}\, \dx \right]
\end{align*}
Since $\frac{p+2}{p+2-q}> \frac{p+2}{2}$ and $R<1$, the previous estimate can be written as follows

\begin{equation}\label{est pulita}
\int_{B_{\frac{R}{2}}}\left|Du(x)\right|^{p+2} \dx\le \frac{c\Arrowvert
u\Arrowvert_{L^\infty\left(B_{R}\right)}^2}{R^\frac{p+2}{p+2-q}}\int_{B_{R}}\left[1+\left|D^2\psi(x)\right|^\frac{p+2}{p+2-q}+
\left|D\psi(x)\right|^\frac{p+2}{p+2-q}+\kappa^{\frac{p+2}{p-q+1}}+\left|Du(x)\right|^p
\right] \dx.
\end{equation}

\noindent Plugging the last inequality in \eqref{stima derivate
seconde} and choosing $\eta\in
C^\infty_0(B_{\frac{R}{2}})$ such that $\eta\equiv 1$ on
$B_{\frac{R}{4}}$ we get

\begin{align*}
\int_{B_{\frac{R}{4}}}&\left|DV_p(Du(x))\right|^2 \dx\le
\frac{c\Arrowvert
u\Arrowvert_{L^\infty\left(B_{R}\right)}^2}{R^\frac{p+2}{p+2-q}}\int_{B_{R}}\left[1+\left|D^2\psi(x)\right|^\frac{p+2}{p+2-q}+
\left|D\psi(x)\right|^\frac{p+2}{p+2-q}+\kappa^{\frac{p+2}{p-q+1}}+\left|Du(x)\right|^p
\right] \dx.
\end{align*}

\noindent that by virtue of estimate \eqref{est boundedness}, gives us the a priori estimate
with

\begin{align}\label{finalestimate}
 \int_{B_{\frac{R}{4}}}\left|DV_p(Du(x))\right|^2 \dx
\le & \frac{c(\Arrowvert \psi\Arrowvert_{L^{\infty}}^2+\Arrowvert
u\Arrowvert_{L^{p^*}\left(B_{R}\right)}^2)}{R^\frac{p+2}{p+2-q}}\notag\\
&\cdot \int_{B_{R}}\left[1+\left|D^2\psi(x)\right|^\frac{p+2}{p+2-q}+
\left|D\psi(x)\right|^\frac{p+2}{p+2-q}+\kappa^{\frac{p+2}{p-q+1}}+\left|Du(x)\right|^p
\right] \dx.
\end{align}

with $c=c(p, q, L, \nu, \mu)$.

\bigskip

{\bf Step 2: The approximation.} \quad Now we conclude the proof by passing to the limit in the approximating problem. The limit procedure is standard { see, e.g., (\cite{CGGP})}.

Let $u \in \mathcal{K}_{\psi}(\Omega)$ be a solution to \eqref{functionalobstacle} and  let ${f_\varepsilon}$ be the sequence obtained applying Lemma \ref{approximation}  to the integrand $f$. Let us fix a ball
$B_R\Subset\Omega$ and let ${u_\varepsilon}\in u+W^{1,p}_0(B_R)$ be the solution to the minimization problem
$$\min\left\{\int_{B_R} {f_\varepsilon}(x,Dv):\,\, v\in \mathcal{K}_{\psi}(B_R)\right\}.$$
By Theorem \ref{thmCGG}, the minimizers ${u_\varepsilon}$ satisfy the  a priori assumptions at \eqref{apriori_estimate_second_derivatives}, i.e. $\left
(\mu^2+\left|Du_\e\right|^2\right)^\frac{p-2}{4}Du_\e\in
W^{1,2}_\mathrm{loc}(\Omega)$, and therefore we are legitimated to use
estimate \eqref{finalestimate} thus obtaing

\begin{align}\label{finalestimateapprox}
 \int_{B_{\frac{R}{4}}}\left|DV_p(Du_\e(x))\right|^2 \dx
\le & \frac{c(\Arrowvert \psi\Arrowvert_{L^{\infty}}^2+\Arrowvert
u_\e\Arrowvert_{L^{p^*}\left(B_{R}\right)}^2)}{R^\frac{p+2}{p+2-q}}\notag\\
&\cdot \int_{B_{R}}\left[1+\left|D^2\psi(x)\right|^\frac{p+2}{p+2-q}+
\left|D\psi(x)\right|^\frac{p+2}{p+2-q}+\kappa^{\frac{p+2}{p-q+1}}+\left|Du_\e(x)\right|^p
\right] \dx.
\end{align}

By the first inequality of growth conditions  at (iii) of Lemma \ref{approximation} and the minimality of $u_\varepsilon$ we get
\begin{eqnarray*}
\int_{B_R} |D{u_\varepsilon}(x)|^p\,\dx &\le& {C(K_0)}\int_{B_R} {f_\varepsilon}(x,D {u_\varepsilon}(x))\, \dx \cr\cr
&\le& {C(K_0)} \int_{B_R} {f_\varepsilon}(x,Du(x))\,\dx\cr\cr
&\le& {C(K_0)}\int_{B_R} f(x,Du(x)) \, \dx \,,
\end{eqnarray*}
where in the last estimate we used the second inequality at (i) of Lemma \ref{approximation}.

Since $f(x, Du)\in L^1_{\mathrm{loc}}(\Omega)$ by assumption,   we deduce, up to subsequences, that there exists $\bar u\in W^{1,p}_0(B_R)+u$ such that
$${u_\varepsilon}\rightharpoonup\bar{u}\quad\quad\text{weakly in }W^{1,p}_0(B_R)+u\,.$$

Note that, since $u_\varepsilon\in \mathcal{K}_\psi$ for every $\varepsilon$ and $\mathcal{K}_\psi$ is a closed set, we have $\bar u\in \mathcal{K}_\psi$. Our next aim is to show that $\bar u$ is a
solution to our obstacle problem over the ball $B_R$.

To this aim fix $\varepsilon_0>0$ and observe that the lower semicontinuity of the functional $w\mapsto \int_{B_R}f_{\varepsilon_0}(x,Dw)\, \dx$, the minimality of $u_\varepsilon$ and the monotonicity of the sequence of
$f_\varepsilon$ yield
\begin{eqnarray*}
\int_{B_R} {f_{\varepsilon_0}}(x,D\bar{u})\, \dx\le \lim_{\varepsilon\to 0}\int_{B_R} f_{\varepsilon_0}(x,Du_{\varepsilon})\, \dx	\le
\int_{B_R} f_{\varepsilon_0}(x,Du)\, \dx\le \int_{B_R} f(x,Du)\, \dx
\end{eqnarray*}
We now use monotone convergence Theorem in the left hand side of previous estimate to deduce that
$$\int_{B_R} f(x,D\bar{u})\,\dx=\lim_{\varepsilon_0\to 0}\int_{B_R} {f_{\varepsilon_0}}(x,D\bar{u})\, \dx\le \int_{B_R} f(x,Du)\, \dx$$
Therefore, we have proved that the limit function $\bar{u}\in W^{1,p}(B_R)+u$ is a solution to the minimization problem
$$\min\left\{\int_\Omega f(x,Dw):w\in W^{1,p}_0(B_R)+u,\,w\in \mathcal{K}_{\psi}\right\}.$$
Since by the strict convexity of the functional the solution is unique, we conclude that $u=\bar u$.
It is quite routine to show that the convergence of ${u_\varepsilon}$ to $u$ is strong in $W^{1,p}_{\mathrm{loc}}(B_R)$.

\noindent The strong convergence of $u_\eps$ to  $u$ in
$W^{1,p}(B_R)$ implies also that $u_\eps$ converges strongly to $u$
in $L^{p^*}(B_R)$ and hence the conclusion follows passing to the limit as $\varepsilon\to 0$ in estimate \eqref{finalestimateapprox}.

\end{proof}

\noindent Acknowledgements: The first and second author have partially been supported by the Gruppo Nazionale per l'Analisi Matematica, la Probabilit\`{a} e le loro Applicazioni (GNAMPA) of the Istituto Nazionale di Alta Matematica (INdAM)

\noindent {\bf A. Gentile}\\
Universit\`{a} degli Studi di Napoli ``Federico II'' \\
Dipartimento di Mat.~e Appl. ``R.~Caccioppoli'',\\
 Via Cintia, 80126 Napoli, Italy

\noindent {\em E-mail address}: andrea.gentile@unina.it

\bigskip
\bigskip

\noindent {\bf R. Giova}\\
\noindent Universit\`{a} degli Studi di Napoli ``Parthenope'' \\
Palazzo Pacanowsky - Via Generale Parisi, 13 \\
 80132 Napoli, Italy

\noindent {\em E-mail address}: raffaella.giova@uniparthenope.it

\bigskip
\bigskip

\noindent {\bf A. Torricelli}\\
\noindent Universit\`{a} degli Studi di Modena e Reggio Emilia\\
\noindent Dipartimento di Scienze Fisiche, Informatiche e Matematiche\\
Via Campi 213/b, 41125 Modena, Italy

\noindent {\em E-mail address}: andrea.torricelli@unipr.it

\bigskip
\bigskip

\end{document}